\documentclass[12pt,a4paper,twoside,fleqn,reqno]{article}
\usepackage{graphicx,enumerate}
\usepackage[dvips]{epsfig}
\usepackage[a4paper,textwidth=4.75in,textheight=7.75in,centering]{geometry}
\usepackage{hyperref} 
\usepackage[T1]{fontenc}
\usepackage[latin1]{inputenc}
\usepackage{amsmath,latexsym,amsfonts,amssymb}
\input epsf

\font\twlmsb=msbm10 at 12pt \font\egtmsb=msbm8 \font\sixmsb=msbm6
\newfam\msbfam
\textfont\msbfam=\twlmsb \scriptfont\msbfam=\egtmsb
\scriptscriptfont\msbfam=\sixmsb
\newtheorem{thm}{Theorem}[section]

\newtheorem{lem}[thm]{Lemma}

\newtheorem{rem}[thm]{Remark}

\title{The Rational Distance Problem for Equilateral Triangles}
\author{Roy Barbara}

\date{}
\begin{document}

\maketitle
\noindent\S-1 \textbf{Abstract}\\
Let \textit{(P)} denote the problem of existence of a point in the plane of a given triangle T, that is at rational distance from all the vertices of T. Answer to \textit{(P)} is positive if T has a rational side and the square of all sides are rational (see [1]). In [2], a complete solution to \textit{(P)} is given for all isosceles triangles with one rational side. In this article, we provide a complete solution to \textit{(P)} for all equilateral triangles.\\
\textit{In all what follows}, $\theta$ denotes an arbitrary positive real number and $T=[\theta]$ denotes the equilateral triangle with side-length $\theta$. For convenience, we say that $\theta$ is "good" (or "suitable") if answer to \textit{(P)} is positive for the triangle $T=[\theta]$. Clearly, the property "$\theta$ is good" is invariant by any rational re-scaling of $\theta$.\\
It turns out that the \textit{good} $\theta$ must have algebraic degree $1,\,2,\,$ or $4$, and they form a subclass of the \textit{positive} bi-quadric numbers, that is, the positive roots of equations of the form\, $x^4+ux^2+v=0$, $\,u,\,v\,\in \mathbb{Q}$. The general form of such numbers is \begin{center}$\sqrt{\alpha\pm\sqrt{\beta}},\;\;\;\;\; \alpha,\,\beta\;\in \mathbb{Q},\;\,\beta\geq \,0,\;\alpha\pm\sqrt{\beta}\,\geq\,0 $\end{center}
that includes positive numbers of the form \begin{center} $\alpha,\,\sqrt{\alpha},\,\alpha\pm\sqrt{\beta}\,\,\sqrt{\alpha}\pm\sqrt{\beta},\;\;\;\;\; \alpha,\,\beta \,\in \mathbb{Q},\;\,\alpha,\,\beta\,\geq \,0.$\end{center}
\textbf{Notations and conventions}: $(x,\,y)$ and $(x,\,y,\,z)$ denote the g.c.d. $\displaystyle \big{(}\frac{x}{p}\big{)}$ denotes legendre's symbol. A triangle with side-lengths $a,\,b,\,c \,$ is denoted by $T=[a,\,b,\,c]$. A triangle is non-degenerated if it has positive area. A radical is non-degenerated if it is irrational.\\ \\
\baselineskip13.5pt
\S-2 \textbf{The results}
\begin{thm}\hspace*{-.3cm}\textbf{0} \,If $\theta$ is good, then, $\theta$ is bi-quadric. More precisely, $\theta^2=\alpha\pm\sqrt{\beta}$ for some $\alpha,\,\beta\in \mathbb{Q},\,\, \beta\geq\,0,$ and $\alpha$ \textit{positive}.
\end{thm}
\begin{thm}\hspace*{-.3cm}\textbf{1} \,Suppose $\theta \notin \mathbb{Q}$ and $\theta^2 \in \mathbb{Q}$. Then, \\
$\theta$ is good $\Leftrightarrow \theta $ has the form $\theta=\lambda\sqrt{p_1...p_r}$ where $\lambda \in \mathbb{Q},\,\lambda\,> \,0,\,r\ge \,1,\, p_1,...,p_r$ are distinct odd primes, $p_i$ is either $3$ or of the form $6k+1$.
\end{thm}
\begin{thm}\hspace*{-.3cm}\textbf{2}
\,Suppose $\theta^2=\alpha\pm \sqrt{\beta},\,\alpha,\,\beta\,\in\mathbb{Q},\,\alpha,\,\beta\,>0,\,\sqrt{\beta}\,\notin \mathbb{Q}$. Then,\\
$\theta$ is good $\;\Leftrightarrow\;$ up to a rational re-scaling of $\theta,\,\,\theta$ is described as follows: \begin{center}$2\theta^2=(a^2+b^2+c^2)\pm\,4\Delta \sqrt{3}$ \end{center}
where $[a,\,b,\,c]$ is a non-degenerated primitive integral triangle with area $\Delta$ such that $4\Delta\,\sqrt{3}\,\notin \mathbb{Q}$.
\end{thm}
\begin{rem}
$\Delta$ is given by Hero's formula,  $\Delta=\sqrt{s(s-a)(s-b)(s-c)},\,s=\frac{1}{2}(a+b+c)$. Equivalently, $4\Delta\sqrt{3}=\sqrt{3(a+b+c)(-a+b+c)(a-b+c)(a+b-c)}$, and the condition $4\Delta\sqrt{3}\,\notin\,\mathbb{Q}$ means that this latter radical is non-degenerated.
\end{rem}
\S-3 \textbf{Proofs of theorems 0 and 1}\\
\textit{Proof of theorem 0:} Suppose $\theta$ good. Let $M$ be a point in the plane of triangle $T=[\theta]$, whose distances from the vertices of $T$ are all rational. the following fundamental relation is well-known (see[3]):
\begin{center} $3(a^4+b^4+c^4+\theta^4)=(a^2+b^2+c^2+\theta^2)^2\;\;\;\;\qquad\quad(\centerdot)$\end{center}
Expanding $(\centerdot)$ yields a relation as $\;\,\theta^4-u\theta^2+v=0$, where $u,\,v\,\in \mathbb{Q}\;$ and $\;u=a^2+b^2+c^2\\>\,0$. Solving for $\theta^2$ yields $\;\,\theta^2=\alpha\pm\sqrt{\beta},$ with $\alpha,\,\beta\,\in \mathbb{Q}\;$ and $\;\alpha=\frac{1}{2}u\,>\,0.$
\begin{lem}\hspace*{-.3cm}\textbf{1}: \,let $q\,>\,1$ be a square-free integer. Then we have:\\
The equation $x^2+3y^2=qz^2$ has a solution in integers  $x,\,y,\,z,$ with $z\ne 0$ if and only if any prime factor of $q$ is either $3$ or of the form $6k+1$
\end{lem}
\textit{\textbf{Proof}}: Suppose first that q has only prime factors as $3$ or $6k+1$. Since the quadratic form $x^2+3y^2,\,x,\,y\,\in \mathbb{Z},$ represents $3$ and every prime $p=6k+1$, and since the set$\{x^2+3y^2,\,x,\,y\,\in \mathbb{Z}\}$ is cosed by multiplication, we conclude that the equation $x^2+3y^2=q.z^2$ has a solution in integers $x,\,y,\,z$ with $z=1.$\\
Conversely, suppose that $x^2+3y^2=q.z^2$ has a solution in integers $x,\,y,\,z,\;\,z\ne 0$. Pick such a solution with $|z|$ minimum. Clearly, $(x,\,y)=1.$ I claim that $q$ is \textit{odd} and has no prime factor $6k-1$. For the purpose of contradiction, we consider two cases:\\
\textbf{case 1}: $q$ is even. Set $q=2w, w$ odd. From $x^2+3y^2=2wz^2$, we see that $x\equiv y (mod. 2).$ As $(x,y)=1,\,x$ and $y$ must be odd, so $x^2+3y^2 \equiv 4 \,(mod. 8)$. Now, $4/2wz^2$ yields $wz^2$ even. But $w$ is odd, hence $z$ is even, so $2wz^2\equiv 0\;(mod. 8)$. We get a contradiction.\\
\textbf{case 2}: $q=p.w$ for some prime $p=6k-1.\,\, x^2+3y^2=pwz^2$ yields $x^2+3y^2\equiv 0 \,(mod. p)$. As $(x,y)=1,\,p$ cannot divide $y$. Hence for some $t\in\mathbb{Z},\,\,yt\equiv 1 \,(mod.p).$ Therefore, $x^2t^2+3y^2t^2\equiv x^2t^2+3\equiv 0 (mod.p),$ so $-3\equiv (xt)^2 \,(mod. p)$. Hence $\displaystyle \big{(}\frac{-3}{p}\big{)}=+1$ contradicting $ p=6k-1$.\\
\begin{lem}\hspace*{-.3cm}\textbf{2}: Let $\theta = \lambda \sqrt{q},\,\lambda\,\in\,\mathbb{Q},\,\lambda > 0,\, q > 1$ square-free integer. We have: \\
$\theta$ is good $\Leftrightarrow$ There are $a,\,b,\,e,\,r,\,s\,\in \mathbb{Q},\,e\ne\,0,$ such that
\begin{center} $a^2+3b^2=q\;\;\;\;\;\quad\quad\qquad\qquad\qquad (1)$ \end{center}
\begin{center} $(a+e)^2+3(b+e)^2=qr^2\;\;\;\;\;\;\qquad(2)$ \end{center}
\begin{center} $(a-e)^2+3(b+e)^2=qs^2\;\;\;\;\;\;\qquad(3)$ \end{center}
\end{lem}
\textit{\textbf{Proof}}: By re-scaling, we take $\theta=2\sqrt{q}$. Let $T=ABC=[\theta]$. Choose a $x-y$ axis to get the coordinates $A(0,\,\sqrt{3q}),\;B(-\sqrt{q},\,0),\;C(\sqrt{q},\,0)$.\\
$\bullet$ Suppose first that $\theta$ is good: There is a point $M=M(x,\,y)$ in the plane of $T$ such that $MA,\, MB,\, MC\,\in \mathbb{Q}$. Clearly, $M\ne A,\,B,\,C$. Set $w=\displaystyle\frac{MA}{q},\,r=\frac{MB}{wq},\,s=\frac{MC}{wq}$. Then, $w,\,r,\,s\,\in \mathbb{Q}-\{0\}.$
The Pythagoras relations are:\\
\begin{center} $\overline{MA}^2=x^2+(y-\sqrt{3q})^2=w^2q^2\;\;\;\;\; (1')$\end{center}
\begin{center} $\overline{MB}^2=(x+\sqrt{q})^2+y^2=w^2q^2r^2\;\;\;\;\; (2')$\end{center}
\begin{center} $\overline{MC}^2=(x-\sqrt{q})^2+y^2=w^2q^2s^2\;\;\;\;\; (3')$\end{center}
Subtracting (2') and (3') yields $\;x=\frac{1}{4}w^2q(r^2-s^2).\sqrt{q},\;$ that is,
\begin{center} $x=\alpha\sqrt{q},\,\;\alpha \in \mathbb{Q}\;\;\;\;\;\;\;\;\;\;\qquad(4)$\end{center}
Then (2') gives $y^2\in\mathbb{Q}$, and then (1') gives $ 2y\sqrt{3q}\in \,\mathbb{Q}$, hence, $y=\gamma \sqrt{3q},\,\gamma\in \mathbb{Q}$.\\
For convenience, we put $\gamma=\beta+1$, obtaining
\begin{center} $y=(\beta+1)\sqrt{3q},\,\,\beta\in\mathbb{Q}\;\;\;\;\quad (5)$\end{center}
Due to (4) and (5), equations (1'), (2'), (3') become after dividing by $q$ :
\begin{center} $\alpha^2+3\beta^2=qw^2$\end{center}
\begin{center} $(\alpha+1)^2+3(\beta+1)^2=qw^2r^2$\end{center}
\begin{center} $(\alpha-1)^2+3(\beta+1)^2=qw^2s^2$\end{center}
Set $a=\displaystyle \frac{\alpha}{w},\,b=\frac{\beta}{w},\,e=\frac{1}{w}$. Dividing by $w^2$, we get precisely relations (1), (2), (3).\\
$\bullet$ Conversely suppose that relations (1), (2), (3) hold with some $a,\,b,\,e,\,r,\,s\,\in \mathbb{Q},\, e\ne 0$. Define point $M=M(x,\,y)$ in the plane of $T$ by \begin{center} $\displaystyle x=\frac{a}{e}\sqrt{q},\;\;\;\; y=(\frac{b}{e}+1)\sqrt{3q}$\end{center}

We may write:\begin{center}$\displaystyle\overline{MA}^2=x^2+(y-\sqrt{3q})^2=q\frac{a^2}{e^2}+3q\frac{b^2}{e^2}=\frac{q}{e^2}(a^2+3b^2)=\frac{q}{e^2}.q=(\frac{q}{e})^2$\end{center}
\begin{center}$\displaystyle\overline{MB}^2=\big{(}(\frac{a+e}{e})\sqrt{q}\big{)}^2+\big{(}(\frac{b+e}{e})\sqrt{3q}\big{)}^2=\frac{q}{e^2}\big{(}(a+e)^2+3(b+e)^2\big{)}=\frac{q}{e^2}.qr^2=(\frac{qr}{e})^2$\end{center}
\begin{center}$\displaystyle\overline{MC}^2=\big{(}(\frac{a-e}{e})\sqrt{q}\big{)}^2+\big{(}(\frac{b+e}{e})\sqrt{3q}\big{)}^2=\frac{q}{e^2}\big{(}(a-e)^2+3(b+e)^2\big{)}=\frac{q}{e^2}.qs^2=(\frac{qs}{e})^2$\end{center}
Therefore, $MA,\,MB,\,MC$ are all rational.\\ \\
\textit{\textbf{Proof of theorem 1}}\\
Let $\theta$ such that $\theta\,\notin \mathbb{Q}$ and $\theta^2\,\in\mathbb{Q}$: $\theta$ can be written as $\theta=\lambda\sqrt{q},\,\lambda\in\,\mathbb{Q},\,\lambda\,>\,0,\,q>1$ square-free integer.\\
$\bullet$ Suppose first that $q$ is even or has a prime factor $6k-1$. By lemma 1, $a^2+3b^2=q,\;\,a,\,b\in \mathbb{Q}$,\; is impossible.\; \textit{Hence}, relation (1) in lemma 2 fails, so $\theta$ is not good.\\
$\bullet$ Suppose now that $q$ has only prime factors as $3$ or $6k+1$. We show that $\theta$ is good using the characterization of lemma 2:\\
By lemma 1, for some $a,\,b\in \mathbb{Q}$, we have $a^2+3b^2=q$. Set $e=\displaystyle -\frac{q}{4b}=\frac{-(a^2+3b^2)}{4b}$,
$r=\frac{a-b}{2b},\;\;\;s=\frac{a+b}{2b}\;\;$. We have\\
$(a+e)^2+3(b+e)^2=(a^2+3b^2)+4e^2+2e(a+3b)=q+\displaystyle\frac{q^2}{4b^2}-\frac{q}{2b}(a+3b)$\\
$\quad=\displaystyle\frac{q}{4b^2}\big{(}4b^2+q-2b(a+3b)\big{)}=\frac{q}{4b^2}(4b^2+a^2+3b^2-2ab-6b^2)$\\
$=\displaystyle\frac{q}{4b^2}(a^2+b^2-2ab)=q\frac{(a-b)^2}{4b^2}=q.r^2 \quad \quad $ and \\
$\displaystyle (a-b)^2+3(b+e)^2=(a^2+3b^2)+4e^2-2e(a-3b)=q+\frac{q^2}{4b^2}+\frac{q}{2b}(a-3b)$\\
$=\displaystyle\frac{q}{4b^2}\big{(}4b^2+q+2b(a-3b)\big{)}=\frac{q}{4b^2}(4b^2+a^2+3b^2+2ab-6b^2)$\\
$=\displaystyle\frac{q}{4b^2}(a^2+b^2+2ab)=q\frac{(a+b)^2}{4b^2}=q.s^2$\\ \\ \\
\newpage
\noindent\S-4 \textbf{\textbf{Proof of theorem 2}}\\
\begin{lem}\hspace*{-.3cm}\textbf{3}:
\,Let $x,\,y,\,z,\,t$ be positive real numbers such that
\begin{center} $3(x^4+y^4+z^4+t^4)=(x^2+y^2+z^2+t^2)^2\;\;\quad\qquad(\circledcirc)$\end{center}
Then, any three of $x,\,y,\,z,\,t$ satisfy the triangle inequality.
\end{lem}
$ $ \\ \textit{\textbf{Proof}}: Since $x,\,y,\,z,\,t$ play symmetric roles, it suffices to show that $x,\,y,\,z$ satisfy the triangle inequality. Write $(\circledcirc$) as \begin{center}$t^4-(x^2+y^2+z^2)t^2+(x^4+y^4+z^4-x^2y^2-y^2z^2-z^2x^2)=0$\end{center}
The discriminant $\bigtriangleup$ of this trinomial in $t^2$ must be non-negative. But, $\bigtriangleup=6(x^2y^2+y^2z^2+z^2x^2)-3(x^4+y^4+z^4)$ that factors as $\;\bigtriangleup =3(x+y+z)(-x+y+z)(x-y+z)(x+y-z)$.\\
Hence, $(-x+y+z)(x-y+z)(x+y-z)\geq 0.$ The reader can easily check (using contraposition) that $x,\,y,\,z$ must satisfy the triangle inequality.\\

\begin{lem}\hspace*{-.3cm}\textbf{4}:\, Let $T=ABC\,=\,\big{[}\theta\big{]}$. Let $a,\,b,\,c$ be positive real numbers satisfying \begin{center}
$3(a^4+b^4+c^4+\theta^4)=(a^2+b^2+c^2+\theta^2)^2$\end{center}
Then, there is a point $M$ in the plane of $T$ such that $\;MA=a,\,MB=b,$ and $MC=c$.
\end{lem}
\textit{\textbf{Proof}}: By lemma 3,\; $a,\,b,$\, and $\theta$ satisfy the triangle inequality. In particular, $a+b\geq\theta$. It folows that the circle $\mathcal{C}(A,a)$ intersects the circle $\mathcal{C}(B,b)$ at two points $M_1$ and $M_2\;\,(M_1=M_2$ if $a+b=\theta)$. Set $c_1=M_1C$ and $c_2=M_2C$. By the fundamental relation $(\centerdot)$ we have $3(a^4+b^4+c_1^4+\theta^4)=(a^2+b^2+c_1^2+\theta^2)^2$ and $3(a^4+b^4+c_2^4+\theta^4)=(a^2+b^2+c_2^2+\theta^2)^2$. Therefore, $c_1^2$ and $c_2^2$ are the roots of the trinomial in $T$
 \begin{center}$ T^2-(a^2+b^2+\theta^2)T+(a^4+b^4+\theta^4-a^2b^2-b^2\theta^2-\theta^2a^2)=0$\end{center}
Since by hypothesis $c^2$ is also a root of this trinomial, we must have $c^2=c_1^2$ or $c^2=c_2^2$. Hence $c=c_1$ or $c=c_2$. Therefore, $a,\,b$ and $c$ are the distances from either point $M_1$ or $M_2\;$ to the vertices $A,\,B$ and $C\;$ of $T$.\\ \\
\textit{\textbf{Proof of theorem 2}}: \\
Let $\theta > 0$ such that $\theta^2=\alpha\pm\sqrt{\beta},\;\,\alpha,\,\beta\,\in\mathbb{Q},\;\alpha,\,\beta\,>0,\;\sqrt{\beta}\notin\mathbb{Q}$.\\
$\bullet$ Suppose first that $\theta$ is good: let $P$ be a point in the plane of $T=ABC=[\theta]$ such that $PA=a,\,PB=b,\,PC=c\;$ are all rational. We have
\begin{center}$3(a^4+b^4+c^4+\theta^4)=(a^2+b^2+c^2+\theta^2)^2\;\;\;\qquad\quad(\centerdot)$\end{center}
By lemma 3, $\;a,\,b,$ and $c\,$ satisfy the triangle inequality. Relation $(\centerdot)$ yields
\begin{center} $\theta^4-U\theta^2+V=0\;$ with $\;U=a^2+b^2+c^2\;$ and $V=a^4+b^4+c^4-a^2b^2-b^2c^2-c^2a^2\;\;(U,V\in\mathbb{Q})$.\end{center}
Solving for $\theta^2$, we get \begin{center}$2\theta^2=(a^2+b^2+c^2)\pm\sqrt{3(a+b+c)(-a+b+c)(a-b+c)(a+b-c)}\;\;\;\quad\quad(\star)$\end{center}
Since $\theta^2$ has algebraic degree $2$, then, the radical in $(\star)$ is non-degenerated. In particular, the triangle $[a,b,c]$ is non degenerated. Select a sufficiently large positive integer $N$ such that $Na,\,Nb,\,Nc$ are all integers and set $D=(Na,\,Nb,\,Nc)$. If we multiply relation $(\star)$ by $\displaystyle\frac{N^2}{D^2}$, this results in replacing in $(\star)\;\theta$ by $\displaystyle\frac{N}{D}.\theta$ and $a,\,b,\,c$ by the integers $\displaystyle\frac{Na}{D},\,\frac{Nb}{D},\,\frac{Nc}{D}$ respectively. As an outcome, we obtain \textit{essentially} the same relation $(\star)$ where $\theta$ has been re-scaled by the rational $\displaystyle\frac{N}{D}$, and where the new symbols $a,\,b,\,c$ represent relatively prime \textit{positive} integers, satisfying the triangle inequality.\\
$\bullet$ Conversely, suppose that for some positive rational $\lambda,\;\theta_0=\lambda.\theta$ is described precisely as in theorem 2. Eliminating the radical
\begin{center}$4\bigtriangleup \sqrt{3}=\sqrt{6(a^2b^2+b^2c^2+c^2a^2)-3(a^4+b^4+c^4)}$\end{center}
in the relation $\;2\theta_0^2=(a^2+b^2+c^2)\pm4\bigtriangleup\sqrt{3}\;\;$ leads to
\begin{center} $3(a^4+b^4+c^4+\theta_0^4)=(a^2+b^2+c^2+\theta_0^2)^2$\end{center}
By lemma 4 there is a point $M$ in the plane of $T=[\theta_0]$ that is at distances $a,\,b,\,c$ from the vertices of $T$. Since $ a,\,b,\,c$ are integers, then, $\theta_0$ is good. Therefore, $\theta=\lambda^{-1}\theta_0$ is also good.\\ \\ \\
We end this article with a few exercises:
\begin{enumerate}
\item Check which are "good" among the radicals: $\sqrt{2},\,\sqrt{3},\,\sqrt{5},\,\sqrt{6},\,\sqrt{7},\,\sqrt{10}$.
\item Show that the positive real number $\theta=\sqrt{25+12\sqrt{3}}$ is "good".
\item Suppose that $2\theta^2=\alpha+\sqrt{\beta},\,\alpha,\,\beta\in\mathbb{Q},\,\alpha,\beta > 0,\sqrt{\beta}\notin\mathbb{Q}$, and $\alpha^2\,<\,\beta$. Show that $\theta$ is not good.
\item Produce solution-points to problem \textit{(P)} for the triangle $T=[\sqrt{3}]$.
\item Let $\theta=\alpha+\beta\sqrt[4]{q}\,>\,0,\;\alpha,\beta\in \mathbb{Q},\,\beta\neq0,\,q\,>\,1$ square-free integer. Show that $\theta$\\ is not good.
\item Suppose that $2\theta^2=\alpha\pm\sqrt{\beta}\,>0,\,\alpha,\,\beta\in\mathbb{Q},\,\alpha,\beta > 0,\sqrt{\beta}\notin\mathbb{Q}$. Write the fraction $\alpha$ in \textit{lowest terms} as $\displaystyle \alpha=\frac{m}{n}\,(m,\,n\;$  positive integers) and suppose that $mn$ has the form $\;mn=4^l(8k+7),\,k,\,l$ non-negative integers. Then, prove that $\theta$ is not good.
\end{enumerate}

\label{lastpage}
\end{document}